\documentclass[12pt,regno]{amsart}
\usepackage{epsf,epsfig,amsmath}
\usepackage{amssymb,latexsym}
\usepackage{tikz}
\usepackage{tikz-cd}
\usepackage{amsthm} 
\usepackage{amsfonts}
\usepackage{graphicx,psfrag}

\numberwithin{equation}{section}

\newtheorem{thm}{Theorem}[section]
\newtheorem{pro}[thm]{Proposition}
\newtheorem{lem}[thm]{Lemma}

\newtheorem{cor}[thm]{Corollary}
\newtheorem*{cor*}{Corollary}

\theoremstyle{remark}\newtheorem{rem}[thm]{Remark}

\theoremstyle{definition}\newtheorem{defi}[thm]{Definition}
\theoremstyle{definition}\newtheorem{exa}[thm]{Example}
\theoremstyle{definition}
\DeclareMathOperator{\Sp}{Sp}
\DeclareMathOperator{\SL}{SL}
\DeclareMathOperator{\SPM}{SPM}

\newcommand{\x}{{\scriptstyle \; \; x \hspace{-6.7pt}\bigcirc \;}}  
\newcommand{\z}{{\scriptstyle \; \; z \hspace{-6.5pt}\bigcirc \;}}  

\title[Pircon kernels and up-down symmetry]{Pircon kernels and up-down symmetry}

\author{Fabrizio Caselli}\author{Mario Marietti}

\address{Fabrizio Caselli, Dipartimento di matematica, Universit\`a di Bologna, Piazza di Porta San Donato 5, 40126 Bologna, Italy}
\address{Mario Marietti, Dipartimento  di Ingegneria Industriale e Scienze Matematiche, Universit\`a Politecnica delle Marche, Via Brecce Bianche, 60131 Ancona,  Italy}

\email{fabrizio.caselli@unibo.it}
\email{m.marietti@univpm.it}


\begin{document}

\begin{abstract} We show that a symmetry property that we call  the up-down symmetry implies that the Kazhdan--Lusztig $R^x$-polynomials of a pircon $P$ are a $P$-kernel, and we show that this property holds in the classical cases.  Then, we enhance and extend to this context a duality of Deodhar in parabolic Kazhdan--Lusztig theory.
\end{abstract}

\maketitle

\section{Introduction}


Kazhdan--Lusztig $R$-polynomials of pircons were introduced in  \cite{Mpirc} in order to provide a  combinatorial generalization not only of  Kazhdan--Lusztig $R$-polynomials of Coxeter groups but also, as instances, of Deodhar's parabolic Kazhdan--Lusztig $R$-polynomials and Kazhdan--Lusztig--Vogan $R$-polynomials and $Q$-polynomials for the action of  $\Sp(2n,\mathbb C)$  on the flag variety of $\SL(2n,\mathbb C)$.

After introducing the Kazhdan--Lusztig $R$-polynomials of a  pircon, the following step is to search for the analog of the Kazhdan--Lusztig $P$-polynomials in this general context. A natural pick is to look at the theory of Kazhdan--Lusztig--Stanley polynomials, introduced by Stanley in \cite{S}.
Indeed, Kazhdan--Lusztig--Stanley polynomials specialize to many interesting objects (see \cite[Sections 6 and 7]{S}). As an example, the  Kazhdan--Lusztig $R$-polynomials of a Coxeter group $W$ form a $W$-kernel whose Kazhdan--Lusztig--Stanley polynomials are the Kazhdan--Lusztig $P$-polynomials of $W$. More generally,  for $H\subseteq S$ and $x\in\{q,-1\}$, the parabolic Kazhdan--Lusztig $R^{H,x}$-polynomials of $W^H$ form a $W^H$-kernel whose Kazhdan--Lusztig--Stanley polynomials are the parabolic Kazhdan--Lusztig $P^{H,x}$-polynomials of $W^H$ (see \cite[Lemma~2.8~(iv) and Proposition~3.1]{Deo87}).

Unfortunately, as shown in \cite[Example~7.2]{Mpirc}, the Kazhdan--Lusztig $R^x$-polynomials of a  pircon $P$ need not to be a $P$-kernel. Therefore, in general, there are no Kazhdan--Lusztig--Stanley polynomials associated with them.

In this work, we provide a sufficient condition for the  Kazhdan--Lusztig $R^x$-polynomials of a  pircon $P$  to be a $P$-kernel, which we call \emph{up-down symmetry}: this roughly says that the $R$-polynomials satisfy a recursive property that is invariant under a flip of $P$ and may have interest in its own right. We show that  Kazhdan--Lusztig--Vogan $R$-polynomials and $Q$-polynomials for the action of  $\Sp(2n,\mathbb C)$  on the flag variety of $\SL(2n,\mathbb C)$ and parabolic Kazhdan--Lusztig $R^{H,x}$-polynomials associated with any Coxeter group satisfy the up-down symmetry. This implies a result by Brenti \cite{BJoA2} on parabolic Kazhdan--Lusztig $R^{H,x}$-polynomials of finite Coxeter groups.

The Hecke algebra $\mathcal H$ of a Coxeter group $W$ also plays a fundamental role in Kazhdan--Lusztig theory. In \cite{Deo87}, Deodhar introduces two $\mathcal H$-modules, one for $x=-1$ and one for $x=q$,  from which he defines the two families of parabolic Kazhdan--Lusztig $R^x$-polynomials and the two families of parabolic Kazhdan--Lusztig $P^x$-polynomials. The aim of a subsequent work, \cite{Deo91}, is to provide a duality  between these two modules and then between the two set ups. Extending Deodhar's duality to the setting of pircons, we introduce two $\mathcal H_P$-modules (where $\mathcal H_P$ is a natural Hecke algebra associated with the pircon $P$) and involutions $\iota^x$ of these modules by means of the up-down symmetry. We then define corresponding Kazhdan--Lusztig bases. We also construct a twisted isomorphism between the two $\mathcal H_P$-modules, which was missing and desired in Deodhar's paper (see Section~\ref{deodhar rivisitato} for further details), and relates the two Kazhdan--Lusztig bases of such modules. We conclude the paper by using these results to find combinatorial recursions to compute the two families of Kazhdan--Lusztig bases and the two families of Kazhdan--Lusztig $P^x$-polynomials.



\section{Notation and preliminaries}
\label{preliminari}

This section reviews the background material that is needed  in the rest of this work.

\subsection{Special partial matchings and pircons.}
Let $P$ be a partially ordered set (poset for short). An element $y\in P$ {\em covers} $x\in P$  if the interval $[x,y]$ coincides with $\{x,y\}$; in this case, we write $x \lhd y$ as well as $y \rhd x$.
If $P$ has a minimum (respectively, a maximum), we denote it by $\hat{0}_P$  (respectively, $\hat{1}_P$). An \emph{order ideal} of $P$ is a subset $I$ of $P$ such that if $y\in I$ and $x\leq y$, then $x\in I$.
The poset $P$ is {\em graded} if $P$ has a minimum and there is a function
$\rho : P \rightarrow {\mathbb N}$ (the {\em rank function}
of $P$) such that $\rho (\hat{0}_P)=0$ and $\rho (y) =\rho (x)
+1$ for all $x,y \in P$ with $x \lhd y$. 
We let $\rho(x,y)=\rho(x)-\rho(y)$, for all $x,y\in P$.
The {\em Hasse diagram} of $P$ is any drawing of the graph having $P$ as vertex set and $ \{ \{ x,y \} \in \binom {P}{2} \colon\,  \text{ either $x \lhd y$ or $y \lhd x$} \}$ as edge set, with the convention that, if $x \lhd y$, then the edge $\{x,y\}$ goes upward from $x$ to $y$. When no confusion arises, we make no distinction between the Hasse diagram and its underlying graph.

%
The following definitions are taken from 
\cite{AH}, and \cite{AHH}, respectively. Given a poset $P$ and $x\in P$, we set $P_{\leq x}=\{y \in P \mid y \leq x\}$.

\begin{defi}\label{accoppiamento parziale}
Let $P$ be a finite poset with $\hat1_P$.
A \emph{special partial matching} of $P$ is an involution $M: P \to P$ such that
\begin{itemize}
  \item $M(\hat1_P) \lhd \hat1_P$,
  \item for all $x\in P$, we have $M(x) \lhd x$, $M(x)=x$, or $M(x)\rhd x$, and
  \item if $x\lhd y$ and $M(x) \neq y$, then $M(x)<M(y)$.
\end{itemize}
\end{defi}
\begin{defi}
A poset $P$ is a \emph{pircon} provided that, for every non-minimal element $x \in P$, the order ideal $P_{\leq x}$ is finite and admits a special partial matching.
\end{defi}
The terminology comes from the fact that a special partial matching without fixed points is precisely a special matching (see \cite{BCM1}, \cite{Mtrans}, and \cite{M} for a definition and for its application to the problem of the combinatorial invariance of  classical and parabolic Kazhdan--Lusztig polynomials) and the fact that pircons relate to special partial matchings in the same way as zircons (see \cite{MJaco}) relate to special matchings. Connected pircons are graded posets (the argument for the zircons in  \cite[Proposition 2.3]{H08} applies also to pircons). 
From now on, when we consider a pircon, we implicitly suppose that it is connected.

\begin{figure}[h]
\begin{center}
$$
\begin{tikzpicture}

\draw[fill=black]{(0,0) circle(2.7pt)};
\draw[fill=black]{(-1.5,1.5) circle(2.7pt)};
\draw[fill=black]{(1.5,1.5) circle(2.7pt)};
\draw[fill=black]{(-2.5,3) circle(2.7pt)};
\draw[fill=black]{(0,3) circle(2.7pt)};
\draw[fill=black]{(2.5,3) circle(2.7pt)};
\draw[fill=black]{(-1.5,4.5) circle(2.7pt)};
\draw[fill=black]{(1.5,4.5) circle(2.7pt)};
\draw[fill=black]{(0,6) circle(2.7pt)};

\draw[dashed, line width=3pt]{(-1.5,1.5)--(0,3)};
\draw[dashed, line width=3pt]{(0,0)--(1.5,1.5)}; 
\draw[dashed, line width=3pt]{(-2.5,3)--(-1.5,4.5)}; 
\draw[dashed, line width=3pt]{(0,6)--(1.5,4.5)}; 

\draw[line width=3pt]{(1.5,1.5)--(2.5,3)};
\draw[line width=3pt]{(-1.5,1.5)--(-2.5,3)};
\draw[line width=3pt]{(0,3)--(1.5,4.5)};
\draw[line width=3pt]{(-1.5,4.5)--(0,6)};

\draw{(1.5,1.5)--(0,3)};
\draw{(0,0)--(-1.5,1.5)}; 
\draw{(2.5,3)--(1.5,4.5)}; 
\draw{(-2.5,3)--(1.5,4.5)}; 
\draw{(0,3)--(-1.5,4.5)};

\draw[line width=2.5pt]{(0,0) circle (.3)};
\draw[dashed, line width=2.5pt]{(2.5,3) circle (.3)};

\end{tikzpicture}$$
\end{center}
\caption{\label{esempio} The thick black special partial matching $M$ and the dashed special partial matching $N$ form a dihedral orbit $\mathcal O$ of rank 3 (top-left) and a chain-like orbit $\mathcal O'$ of rank 2 (bottom-right). We have $m(\mathcal O)=m(\mathcal O')=3$.} 
\end{figure}
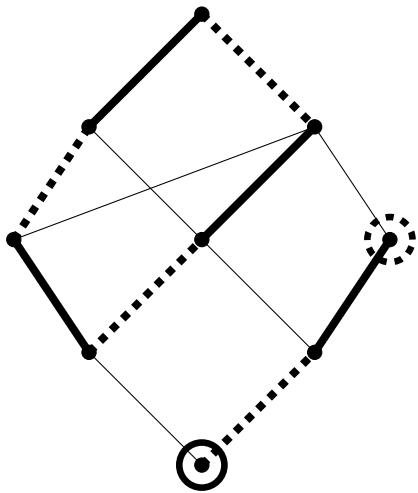

Given a poset $P$ and $w\in P$, we say that $M$ is a matching of $w$ if $M$ is a matching of $P_{\leq w}$, and we denote by $\SPM_w$ the set of all special partial matchings of $w$. Hence, if $P$ is a pircon then $\SPM_w\neq \emptyset$ for all $w\in P\setminus \{\hat{0}_P\} $. 
In pictures, we visualize a special partial matching $M$ of a poset $P$ by taking the Hasse diagram of $P$ and coloring in the same way, for all $x\in P$, either the edge $\{x,M(x)\}$ if $M(x)\neq x$, or a circle around $x$ if $M(x)=x$.
Figure~\ref{esempio} shows two special partial matchings of a poset, one colored with black thick lines which fixes only the bottom element $\hat{0}_P$, and one colored with dashed lines which also fixes only one element. 

\subsection{Coxeter groups}
\label{gruppi di coxeter}
We  follow   \cite{BB} for undefined  terminology concerning Coxeter groups.
We fix our notation on a Coxeter system $(W,S)$ in the following list:
\smallskip 

{\renewcommand{\arraystretch}{1.2}
$
\begin{array}{@{\hskip-1.3pt}l@{\qquad}l}
m(s,t) &  \textrm{the entry of the Coxeter matrix of  $(W,S)$ in position $(s,t)\in S\times S$}, 
\\
e &  \textrm{identity of $W$}, 
\\
\ell  &  \textrm{the length function of $(W,S)$},
\\
D_R(w) & =\{ s \in S : \; \ell(w  s) < \ell(w ) \},  \textrm{ the right descent set of $w\in W$},
\\
D_L(w) & =\{ s \in S : \; \ell( sw) < \ell(w ) \},  \textrm{ the left descent set of $w\in W$},
\\
W_H & \textrm{ the parabolic subgroup of $W$ generated by $H\subseteq S$},
\\
W^H &=\{ w \in W \, : \; D_{R}(w)\subseteq S\setminus H \}, \textrm{ the set of minimal left coset representatives},
\\
\leq & \textrm{ Bruhat order on $W$ (as well as any other order on a poset $P$)},
\\
\textrm{$[u,v]$} & =\{ w \in W \, : \; u \leq w \leq v \}, \textrm{ the (Bruhat) interval generated by $u,v\in W$},
\\
\textrm{$[u,v]^{H}$} &  = \{ z \in W^{H}: \;  u \leq z \leq v \}, \textrm{ the parabolic (Bruhat)  interval generated by $u,v\in W^H$}.

\end{array}$
\bigskip

The \emph{Hecke algebra} of W, denoted $\mathcal H(W)$,  is the $\mathbb Z[q^{\frac{1}{2}},q^{-\frac{1}{2}}]$-algebra generated by $\{T_s : s \in S\}$ subject to the braid relations 
$$\underbrace{\cdots T_sT_rT_s}_{\text{$m(s,r)$ terms}}= \underbrace{\cdots T_rT_sT_r}_{\text{$m(s,r)$ terms}}\quad \text{ for all $s,r \in S$} $$
and the quadratic relations
$$T^2_s=(q-1)T_s+q \quad \text{ for all $s \in S$}.$$
For $w\in W$, denote by $T_w$ the product $T_{s_1}T_{s_2}\cdots T_ {s_k}$, where $s_1 s_2\cdots s_k$ is a reduced expression for $w$. The element $T_w$ is independent from the chosen reduced expression. The Hecke algebra $\mathcal H(W)$ is a free $\mathbb Z[q^{\frac{1}{2}},q^{-\frac{1}{2}}]$-module having the set $\{T_w : w \in W\}$ as a basis and multiplication uniquely determined by
$$ T_sT_w =\left\{
\begin{array}{ll}
T_{sw}, & \text{ if $sw >w$,} \\
qT_{sw} + (q-1)T_w, & \text{  if  $sw <w$,}
\end{array}
\right.
$$
for all $w \in W$ and $s\in S$. 

For any $w\in W$, the element $T_w$ is invertible; for example, if $s\in S$ then $T_{s}^{-1}= q^{-1} T_s + (q^{-1}-1)$. We denote by $\iota$ the involution defined by $\iota (\sum a_w \; T_w)=\sum \overline{a}_w \; T_{w^{-1}}^{-1}$, where $a\mapsto \overline{a}$ is the involution of the ring $\mathbb Z[q^{\frac{1}{2}},q^{-\frac{1}{2}}]$ sending  $q^{\frac{1}{2}}$ to $q^{-\frac{1}{2}}$.
\bigskip

\subsection{Kazhdan--Lusztig $R$-polynomials for pircons}
\label{orbits}
In this subsection, we recall the results from \cite{Mpirc} that are needed later. Some of them are valid (with the same proof) not only for the special partial matchings but also for the larger class of \emph{quasi special partial matchings}, that we introduce here. Since we actually need these results for this larger class, we state them directly in this more generalized form. 
\begin{defi}\label{accoppiamento quasi parziale}
Let $P$ be a  poset.
A \emph{quasi special partial matching} of $P$ is an involution $M: P \to P$ such that
\begin{itemize}
  \item for all $x\in P$, we have $M(x) \lhd x$, $M(x)=x$, or $M(x)\rhd x$, and
  \item if $x\lhd y$ and $M(x) \neq y$, then $M(x)<M(y)$.
\end{itemize}
\end{defi}

So, a quasi special partial matching $M$ of a finite poset $P$ with $\hat 1_P$ is a special partial matching if and only if  $M(\hat1_P) \lhd \hat1_P$. Notice that in the definition of a quasi special partial matching the poset $P$ needs not be finite.

For the proof of the following result (in the case of special partial matchings), see \cite[Lemma~5.2]{AHH}.
\begin{lem}[Lifting property for quasi special partial matchings]
\label{lifting}
Let  $M$ be a quasi special partial matching of a poset $P$. If $x,y \in P$ with $x<y$ and $M(y) \leq y$, then
\begin{itemize}
  \item[(i)]   $M(x) \leq y$,
  \item[(ii)]  $M(x) \leq x \implies M(x)<M(y)$, and
  \item[(iii)] $M(x) \geq x \implies x \leq M(y)$.
\end{itemize}
\end{lem}

Lemma~\ref{lifting} implies the following.
\begin{lem}
\label{restringe}
Let $M$ be a quasi special partial matching of a poset $P$, and $x,y \in P$, with $x\leq y$.  If $M(y)\leq y$ and $M(x) \geq x$, then $M$ restricts to a quasi special partial matching of the interval $[x,y]$.
\end{lem}

We say that an interval $[u, v]$ in a poset $P$ is {\em dihedral} if it is isomorphic to an interval in a Coxeter group with  two Coxeter generators, ordered by Bruhat order (see Figure~\ref{diheinte}).
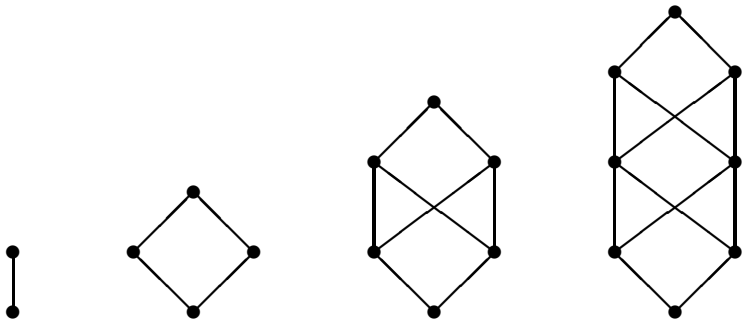
\begin{figure}[h]
 \setlength{\unitlength}{8mm}
\begin{center}
\begin{picture}(12,6)
\thicklines
\put(0,0){\line(0,1){1}}
\put(0,0){\circle*{0.2}}
\put(0,1){\circle*{0.2}}

\put(3,0){\line(1,1){1}}
\put(3,0){\line(-1,1){1}}
\put(4,1){\line(-1,1){1}}
\put(2,1){\line(1,1){1}}
\put(3,0){\circle*{0.2}}
\put(3,2){\circle*{0.2}}
\put(2,1){\circle*{0.2}}
\put(4,1){\circle*{0.2}}

\put(7,0){\line(1,1){1}}
\put(7,0){\line(-1,1){1}}
\put(8,1){\line(0,1){1.5}}
\put(6,1){\line(0,1){1.5}}
\put(6,2.5){\line(1,1){1}}
\put(8,2.5){\line(-1,1){1}}
\put(6,1){\line(4,3){2}}
\put(8,1){\line(-4,3){2}}
\put(7,0){\circle*{0.2}}
\put(8,1){\circle*{0.2}}
\put(6,1){\circle*{0.2}}
\put(6,2.5){\circle*{0.2}}
\put(8,2.5){\circle*{0.2}}
\put(7,3.5){\circle*{0.2}}

\put(11,0){\line(1,1){1}}
\put(11,0){\line(-1,1){1}}
\put(12,1){\line(0,1){3}}
\put(10,1){\line(0,1){3}}
\put(10,1){\line(4,3){2}}
\put(12,1){\line(-4,3){2}}
\put(10,2.5){\line(4,3){2}}
\put(12,2.5){\line(-4,3){2}}
\put(10,4){\line(1,1){1}}
\put(12,4){\line(-1,1){1}}
\put(11,0){\circle*{0.2}}
\put(12,1){\circle*{0.2}}
\put(10,1){\circle*{0.2}}
\put(12,2.5){\circle*{0.2}}
\put(10,2.5){\circle*{0.2}}
\put(10,4){\circle*{0.2}}
\put(12,4){\circle*{0.2}}
\put(11,5){\circle*{0.2}}
\end{picture}
\end{center}
\caption{\label{diheinte} Dihedral intervals of rank 1,2,3,4}
\end{figure}

Given a pircon $P$ and two quasi special partial matchings  $M$ and $N$ of $P$,  we denote by $\langle M,N \rangle $ the group of permutations of $P$ generated by $M$ and $N$. Furthermore, we denote by $\langle M,N \rangle (u)$ the orbit of an element $u\in P$ under the action of $\langle M,N \rangle $.

We say that an orbit $\mathcal O$ of the action of $\langle M,N \rangle $ is 
\begin{itemize}
\item \emph{dihedral}, if  $\mathcal O$ is isomorphic to a dihedral interval and $u\notin \{M(u),N(u)\}$ for all $u\in \mathcal O$ (see Figure~\ref{esempio}, top-left),
\item \emph{chain-like}, if  $\mathcal O$ is isomorphic to a chain, $\hat{0}_{\mathcal O} \in \{M(\hat{0}_{\mathcal O}),N(\hat{0}_{\mathcal O})\}$, and $\hat{1}_{\mathcal O} \in \{M(\hat{1}_{\mathcal O}),N(\hat{1}_{\mathcal O})\}$ (see Figure~\ref{esempio}, bottom-right).
\end{itemize}

Note that an orbit with two elements $w$ and  $N(w)=M(w)\neq w$ is dihedral, whereas an orbit with two elements $w=N(w)$ and  $M(w)=NM(w)\neq w$ is chain-like (see Figure~\ref{rango1}). This is the only case when a dihedral orbit and a chain-like orbit are isomorphic as posets.

\begin{figure}[h]
\begin{center}
$$
\begin{tikzpicture}
 \draw[fill=black]{(2,2) circle(3.5pt)};
\draw[fill=black]{(2,0) circle(3.5pt)};
 \draw[line width=2.5pt]{(1.95,0)--(1.95,2)}; 
  \draw[line width=2.5pt, dashed]{(2.05,0)--(2.05,2)}; 
 \draw[fill=black]{(5,2) circle(3pt)};
\draw[fill=black]{(5,0) circle(3pt)};
\draw[dashed, line width=2pt]{(5,0) circle (.25)};
 \draw[line width=3pt]{(5,0)--(5,2)}; 
\draw[line width=2pt, dashed]{(5,2) circle (.25)};
\end{tikzpicture}$$
\end{center}
\caption{\label{rango1} A dihedral orbit of rank 1 (left)  and a chain-like orbit of rank 1 (right).}
\end{figure}
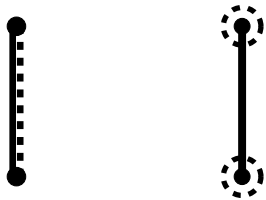

\begin{lem}
\label{noorbita}
Fix a pircon $P$.  Let $M$ and $N$ be two quasi special partial matchings of an element $w$ in $P$. 
Every orbit $\mathcal O$ of $\langle M,N \rangle $ is either dihedral or chain-like. Moreover, $\mathcal O$  is an interval in $P$, i.e. $\mathcal O=[\hat{0}_{\mathcal O},\hat{1}_{\mathcal O}]$.
\end{lem}

Recall that, given a pircon $P$ and $w\in P$, we denote by $\SPM_w$ the set of all special partial matchings of $w$. Notice that, if $M,N \in \SPM_w$, then the orbit $\langle M,N \rangle (w)$ is dihedral by the definition of a special partial matching and Lemma~\ref{noorbita}. The following definition is a generalization of \cite[Definition~3.1]{BCM2}.
Given  two quasi special partial matchings $M$ and $N$ and a finite orbit $\mathcal O$ of $\langle M,N \rangle $, we let (see Figure~\ref{esempio})
$$m(\mathcal O)= \left\{ \begin{array}{ll}
\mbox{the rank of $\mathcal O$}, & \mbox{if $\mathcal O$ is dihedral,} \\
\mbox{the rank of $\mathcal O + 1$}, & \mbox{if $\mathcal O$ is  chain-like.}
\end{array} \right. $$

\begin{defi}
Let $P$ be a pircon and  $w \in P$. We say that  two quasi special partial matchings $M,N$ of $w$ are
\emph{strictly coherent} provided that  $$m(\mathcal O)\textrm{ is a divisor of }m(\langle M,N \rangle (w))$$ 
for every orbit $\mathcal O$  of $\langle M,N \rangle $. Moreover, we say that two special partial matchings $M,N$ of $w$ are \emph{coherent} provided there exists a sequence $M_0,M_1, . . . , M_k$ of special partial
matchings of $w$ such that $M_0 = M$, $M_k = N$, and $M_i$ and $M_{i+1}$ are strictly coherent for all $i = 0,1, . . . , k-1$. 
\end{defi}

Notice that, while the definition of coherent is given for special partial matchings, the definition 
  of strictly coherent is valid more generally  for quasi special partial matchings.

By the definition of a pircon,  we can fix one special partial matching of $v$  for each $ v\in P\setminus \{\hat{0}_P\} $. If  $\mathcal M$ is the set of such fixed special partial matchings, then we call the pair $(P,\mathcal M)$ a \emph{refined pircon} and $\mathcal M$ a \emph{refinement} of $P$.

\begin{defi}
\label{klmm}
Let $x\in \{q,-1\}$.
Let $(P,\mathcal M)$ be a refined pircon, where $\mathcal M=  \{M_v \in \SPM_v: v\in P\setminus \{\hat{0}\} \}$. The family of \emph{Kazhdan--Lusztig $R^{x}$-polynomials} $\{R^{x}_{u,w}{}\}_{u,w\in P}\subseteq \mathbb Z[q]$ of $(P,\mathcal M)$ (or $R^{x}$-polynomials for short) is the unique family of polynomials satisfying  the following properties:
\begin{itemize}
\item if $u\not\leq w$ then $R_{u,w}^{x} (q)=0$,
\item $R_{w,w}^{x} (q)=1$ for all $w\in P$,
\item if $u\leq w$ then 
\begin{equation}
\label{klperpirconi}
 R_{u,w}^{x} (q)= \left\{ \begin{array}{ll}
R_{M_w(u),M_w(w)}^{x}(q), & \mbox{if $M_w(u)  \lhd u$,} \\
(q-1)R_{u,M_w(w)}^{x}(q)+qR_{M_w(u),M_w(w)}^{x}(q), & \mbox{if $M_w(u) \rhd u$,} \\
(q-1-x)R_{u,M_w(w)}^{x}(q), & \mbox{if $M_w(u) = u$.} 
\end{array} \right. 
\end{equation}
\end{itemize}
\end{defi} 
\begin{exa}\label{R-chain}
	If $P$ is a chain and $\mathcal M$ is any refinement of $P$, then the $R^x$-polynomials of $(P,\mathcal M)$ are
	\[
	R^x_{u,v}(q)=(q-1)(q-1-x)^{\rho(u,v)-1}
	\] 
	for all $u,v\in P$ with $u<v$.
	\end{exa}
	
In general, the Kazhdan--Lusztig $R^{x}$-polynomials of a refined pircon $(P,\mathcal M)$  depend on the refinement $\mathcal M$ (see \cite[Remark~5.3]{Mpirc}). 

The two families of $R$-polynomials satisfy the following properties.
\begin{pro}
\label{relazione}
Let $(P,\mathcal M)$ be a refined pircon with rank function $\rho$. If $u,w\in P$ with $u\leq w$, then 
\begin{enumerate}
\item $\deg R^{-1}_{u,w}(q)=\rho(u,w)$,
\item $R^{q}_{u,w}(0)=(-1)^{\rho(u,w)}$,
\item $R^{x}_{u,w}(q)=(-q)^{\rho(u,w)} \; R^{z}_{u,w}(q^{-1})$, where $\{x,z\}=\{q,-1\}$
\end{enumerate}
\end{pro}

Let $(P,\mathcal M )$ be a refined pircon and $w\in P$. Mimicking \cite[Subsection~4.1]{Mtrans}, we say that a special partial matching $M$ of $w$ \emph{calculates} the Kazhdan--Lusztig $R^{x}$-polynomials of $(P,\mathcal M)$ (or  is \emph{calculating}, for short) provided that, for all $u \in P$, $u \leq w$, the following holds:
\begin{equation}\label{calculat}
 R_{u,w}^{x} (q)= \left\{ \begin{array}{ll}
R_{M(u),M(w)}^{x}(q), & \mbox{if $M(u)  \lhd u$,} \\
(q-1)R_{u,M(w)}^{x}(q)+qR_{M(u),M(w)}^{x}(q), & \mbox{if $M(u) \rhd u$,} \\
(q-1-x)R_{u,M(w)}^{x}(q), & \mbox{if $M(u) = u$.} 
\end{array} \right. 
\end{equation}
Thus the matchings of $\mathcal M$ are calculating by definition.

\begin{defi}
Let $(P,\mathcal M )$ be a refined pircon. We say that a quasi special partial matching $M$ of an order ideal $I$ of $P$ is \emph{strongly calculating} provided that the restriction of $M$ to $P_{\leq z}$ is calculating for all $z\in I$ such that $M(z)\lhd z$.
\end{defi}
Notice that Eq.~\eqref{calculat} holds also if $u\not \leq w$, since in this case all terms are zero by Lemma~\ref{lifting} and Definition~\ref{klmm}.
Notice also that, by Lemma~\ref{restringe}, the restriction of a quasi special partial matching $M$ to $P_{\leq z}$ is indeed a special partial matching for all $z\in P$ such that $M(z)\lhd z$.

The following definition slightly differs from \cite[Definition~5.6]{Mpirc} but is more suitable for our purposes.
\begin{defi}
\label{sistema}
We say that $(P,S)$ is a \emph{pircon system} provided that
\begin{enumerate}
\item $P$ is a pircon,
\item $S$ is a set of quasi special partial matchings of order ideals of $P$, 
\item
\label{tre}  for all $w\in P\setminus \{\hat{0}_P\}$, there exists $M\in S$ such that $M(w)$ is defined and $M(w)\lhd w$,
\item for all $w\in P$ and all $M,N\in S$ such that $M(w)$ and $N(w)$ are defined and satisfy $M(w)\lhd w$ and $N(w)\lhd w$, the restrictions of $M$ and $N$ to $P_{\leq w}$ are coherent.
\end{enumerate}
\end{defi}

Although our definition of pircon system is different from that of \cite{Mpirc}, the proof of the following result is equal to the proof of \cite[Corollary~5.7]{Mpirc}.

\begin{thm}
\label{unici}
Let $(P,S)$ be a pircon system and $x\in\{q,-1\}$. All refinements $\mathcal M$ of $P$, with $\mathcal M \subseteq S$, yields the same family of Kazhdan--Lusztig $R^x$-polynomials (for which, all matchings in $S$ are strongly calculating).
\end{thm}

\begin{defi}
\label{dircone}
A pircon $D$ is a \emph{dircon} provided that any two special partial matchings $M,N \in \SPM_w$ are coherent,  for all $w \in D$. 
\end{defi}
In other words, a pircon $D$ is a dircon if and only if $(D, \bigcup_{w\in P\setminus \{\hat{0}_P\}}\SPM_w)$ is a  pircon system. By Theorem~\ref{unici}, for both $x=q$ and $x=-1$, a dircon has a unique family of Kazhdan--Lusztig $R^x$-polynomials. The terminology comes from the fact that dircons relate to pircons in the same way as diamonds relate to zircons (see \cite[Definition~3.2]{BCM2}).

\subsection{Stanley's kernels}
\label{KLS}
 Let us briefly recall from \cite{S} the definition of $P$-kernel and Kazhdan--Lusztig--Stanley polynomials.

Let $P$ be a locally finite graded poset, with rank function $\rho$. The \emph{incidence algebra} of $P$ over the polynomial ring  $\mathbb R[q]$, denoted
$I(P)$, is the associative algebra of functions $f$ assigning to each nonempty interval $[u, v]$ an element  $f_{u,v}(q)\in \mathbb R[q]$ (denoted also simply by $f_{u,v}$ when no confusion arises) with usual sum and convolution product: $(f+g)_{u,v}=f_{u,v}+g_{u,v}$ and $(f \cdot g)_{u,v}=\sum_{z:u\leq z\leq v} f_{u,z} \;g_{z,v}$, for all $f,g \in  I(P)$ and all $u,v\in P$ with $u\leq v$. The  identity element of  $I(P)$ is the \emph{delta function} $\delta$, defined by $\delta_{u,v}=\left\{\begin{array}{ll}
1 & \text{if $u=v$,} \\
0 & \text{if $u<v$.}
\end{array} \right.$

If $f\in I(P)$, we also let $f_{u,v}=0$ whenever $u\not\leq v$, and so  $(f \cdot g)_{u,v}=\sum_{z\in P} f_{u,z} \;g_{z,v}$, for all $f,g \in  I(P)$. An element $f\in I(P)$ is invertible if and only if $f_{u,u} \in \mathbb R \setminus \{0\}$ for all $u\in P$.  We say that $f\in I(P)$ is unitary if $f_{u,u}=1$, for all $u\in P$.
Let 
\begin{itemize}
\item $I'(P)=\{f \in I(P): \deg f_{u,v} \leq \rho(u,v)), \text{ for all $u,v\in P$ with $u\leq v$}\}$, 
\item  $I_{\frac{1}{2}}(P)=\{f\in I'(P) \text{ unitary}: \deg f_{u,v}<\frac{1}{2}\rho(u,v), \text{ for all  $u,v\in P$ with $u< v$}\}$.
\end{itemize}
Note that $I'(P)$ is a subalgebra of $I(P)$, closed under taking inverse. 
Given $f \in I'(P)$,  we denote by $\widetilde{f}$ the element of $I'(P)$  such that $\widetilde{f}_{u,v}(q)=q^{\rho(u,v)} f_{u,v}(q^{-1})$, for all  $u,v\in P$ with $u\leq v$. 
Notice that the map $\; \widetilde{\text{}} \;$ is an involution on $I'(P)$.

A unitary element $K\in I(P)$ is a \emph{$P$-kernel} if there exists
an invertible element $f \in I(P)$ such that $K \cdot f= \widetilde{f}$.
Such an  element $f\in I(P)$  is called \emph{invertible $K$-totally acceptable} function in \cite{S}.  See \cite[Theorem~6.5, Proposition~6.3, Corollary~6.7]{S} for a proof of the next result.
\begin{thm}\label{Pkernelequiv}
Let $P$ be a locally finite graded poset.
\begin{enumerate}
\item A unitary $K\in I'(P)$  is a $P$-kernel if and only if $K \cdot \widetilde{K}=\delta$.
\item There is a bijection from the set of $P$-kernels of $I'(P)$ to $ I_{\frac{1}{2}}(P)$ that assigns to $K$ an invertible $K$-totally acceptable function.
\end{enumerate} 
\end{thm} 

Let $K\in I'(P)$ be a $P$-kernel. Following \cite{BJoA}, we refer to the unique invertible $K$-totally acceptable function of $ I_{\frac{1}{2}}(P)$ as the \emph{Kazhdan--Lusztig--Stanley polynomials} of $K$.

\section{Up-down symmetry and $P$-polynomials of pircons} \label{sec:4}

In this section, we introduce and study the concept of \emph{up-down symmetry}. In particular, we relate it to Stanley's kernels and to the existence of the $P$-polynomials of a pircon. Notice that, in general, the family $\{R^x_{u,v}(q)\}_{u,v\in P}$ of  $R^x$-polynomials needs not be a $P$-kernel (see \cite[Example~7.2]{Mpirc}).

Recall that we write  $\overline{q}$ for $q^{-1}$ and, for a polynomial $f$, we write $\overline{f}$ for $f(\overline{q})$. Recall also the three instances of the recursive formula for a strongly calculating quasi special partial matching $M$: for all $u$ and $w$ such that $M(u)$ and $M(w)$ are defined and $M(w)\lhd w$,  we have
\begin{align}\label{RecR}
R_{u,w}^{x} (q)= \begin{cases} 
R_{M(u),M(w)}^{x}(q), & \mbox{if $M(u)  \lhd u$, \hspace{1cm}(a)}\\
(q-1)R_{u,M(w)}^{x}(q)+qR_{M(u),M(w)}^{x}(q), & \mbox{if $M(u) \rhd u$, \hspace{1cm}(b)}\\
(q-1-x)R_{u,M(w)}^{x}(q), & \mbox{if $M(u) = u$. \hspace{1cm}(c)}
\end{cases} 
\end{align}
The next definition involves analogous formulas obtained by exchanging the roles of the elements $u$ and $w$.
\begin{defi}Let $x\in \{q,-1\}$. The family of Kazhdan--Lusztig $R^x$-polynomials of a pircon system $(P,S)$  satisfies the up-down symmetry if, for all $M\in S$ and for all $u,w\in P$ such that $M(u)$ and $M(w)$ are defined and $M(u)\rhd u$, we have
\begin{align}\label{RecR'}
R_{u,w}^{x} (q)= \begin{cases} 
R_{M(u),M(w)}^{x}(q), & \mbox{if $M(w)  \rhd w$, \hspace{1cm}(a$'$)}\\
(q-1)R_{M(u),w}^{x}(q)+qR_{M(u),M(w)}^{x}(q), & \mbox{if $M(w) \lhd w$, \hspace{1cm}(b$'$)}\\
(q-1-x)R_{M(u),w}^{x}(q), & \mbox{if $M(w) = w$. \hspace{1cm}(c$'$)}
\end{cases} 
\end{align}	
	\end{defi}
\begin{lem}
\label{cesena}
	Let $x\in \{q,-1\}$. The family of Kazhdan--Lusztig $R^x$-polynomials of a pircon system $(P,S)$  satisfies the up-down symmetry if and only if, for all $M\in S$ and all $u,w\in P$ such that $M(u)\rhd u$ and $M(w)=w$, we have
	\begin{equation}
	\label{equazione}
	R_{u,w}^{x}=(q-1-x) \, R_{M(u),w}^{x}.
	\end{equation}
	\end{lem}
\begin{proof}
	It is a trivial check that the two properties (a$'$) and (b$'$)  in Eq.~\eqref{RecR'} follow by the two properties (a) and (b) in Eq.~\eqref{RecR}. The result follows.
\end{proof}	


\begin{thm}\label{arePkernel}
Let $(P,S)$ be a pircon system and $x\in\{q,-1\}$. If the family $\{R^x_{u,v}(q)\}_{u,v\in P}$ of  $R^x$-polynomials satisfies the up-down symmetry, then it defines a $P$-kernel.
\end{thm}
\begin{proof}By Theorem~\ref{Pkernelequiv}, we have to show that
 \[
  \sum_{z\in P}R^x_{u,z} \; q^{\rho(z,v)} \; \overline{R^x_{z,v}}=0
 \]
holds  for all $u,v\in P$ with $u< v$. We proceed by induction on $\rho(v)$. The case $\rho(v)=1$ is trivial. Suppose $\rho(v)>1$. Moreover, we may suppose also $\rho(u,v)>1$ since the case $\rho(u,v)=1$ is trivial as well. We fix a quasi special partial matching $M\in S$ such that $M(v)$ is defined and $M(v)\lhd v$.
 For notational convenience, for all $u$ and $v$, we let
 \begin{itemize}
  \item $\Delta_{=}(u,v)=\displaystyle\sum_{z:\, M(z)= z}R^x_{u,z} \; q^{\rho(z,v)} \; \overline{R^x_{z,v}}$,
  \item $\Delta_{\lhd}(u,v)=\displaystyle\sum_{z:\, M(z)\lhd z}R^x_{u,z} \; q^{\rho(z,v)} \; \overline{R^x_{z,v}}$,
  \item $\Delta_{\rhd}(u,v)=\displaystyle\sum_{z:\, M(z)\rhd z}R^x_{u,z} \; q^{\rho(z,v)} \; \overline{R^x_{z,v}}$,
  \item  $\Delta(u,v)=\Delta_{=}(u,v)+\Delta_{\lhd}(u,v)+\Delta_{\rhd}(u,v)$.
 \end{itemize}
   We need to prove that $\Delta(u,v)=0$, for all $u,v\in P$ with $u< v$. 
   
%
 We split the proof into three cases depending on whether $M(u)\lhd u$, $M(u)\rhd u$, or $M(u)=u$.

 Suppose $M(u)\lhd u$. We have
 
 \begin{align}
\label{onne} \Delta_{\lhd}(u,v)&=\sum_{z:\, M(z)\lhd z}R^x_{M(u),M(z)} \; q^{\rho(z,v)} \; \overline{R^x_{M(z),M(v)}}=\sum_{w:\, M(w)\rhd w}R^x_{M(u),w} \; q^{\rho(M(w),v)} \; \overline{R^x_{w,M(v)}}\\
 \nonumber &= \sum_{w:\, M(w)\rhd w}R^x_{M(u),w} \; q^{\rho(w,M(v))} \; \overline{R^x_{w,M(v)}}= \Delta_{\rhd}(M(u),M(v))\\
\nonumber  &=- \Delta_{\lhd}(M(u),M(v))- \Delta_{=}(M(u),M(v))
 \end{align}
 using the induction hypothesis on the interval $[M(u),M(v)]$.
 We also have 
 \begin{align}
 \label{twwo} \Delta_{\lhd}(M(u),M(v))&= \sum_{z:\, M(z)\lhd z}R^x_{M(u),z} \; q^{\rho(z,M(v))} \; \overline{R^x_{z,M(v)}}\\
\nonumber  &= \sum_{z:\, M(z)\lhd z}\big((q-1)R^x_{M(u),M(z)}+qR^x_{u,M(z)}\big) \; q^{\rho(z,M(v))} \; \overline{R^x_{z,M(v)}}\\
 \nonumber  &=\sum_{z:\, M(z)\lhd z}(q-1)R^x_{u,z} \; q^{\rho(z,M(v))} \; \overline{R^x_{z,M(v)}}+  \sum_{w:\, M(w)\rhd w}qR^x_{u,w} \; q^{\rho(M(w),M(v))} \; \overline{R^x_{M(w),M(v)}}\\
 \nonumber  &=(q-1)\Delta_{\lhd}(u,M(v))+q\sum_{w:\, M(w)\rhd w}R^x_{u,w}\;  q^{\rho(M(w),M(v))} \;  \overline{R^x_{M(w),M(v)}}.
 \end{align}
 
On the other hand
\begin{align}
\label{thrree} \Delta_{\rhd}(u,v)&=\sum_{z:\, M(z)\rhd z}R^x_{u,z} \;  q^{\rho(z,v)} \;  \overline{R^x_{z,v}}\\
\nonumber  &= \sum_{z:\, M(z)\rhd z}R^x_{u,z} \;  q^{\rho(z,v)}\big( \overline{( q-1) \; R^x_{z,M(v)}}+ \overline{qR^x_{M(z),M(v)}} \big)\\
\nonumber  &=\sum_{z:\, M(z)\rhd z}(1-q)R^x_{u,z}\;  q^{\rho(z,M(v))}\;  \overline{R^x_{z,M(v)}}+\sum_{z:\, M(z)\rhd z}qR^x_{u,z}\;  q^{\rho(M(z),M(v))}\;  \overline{R^x_{M(z),M(v)}}\\
\nonumber  &=(1-q) \Delta_{\rhd}(u,M(v))+q\sum_{z:\, M(z)\rhd z}R^x_{u,z}\;  q^{\rho(M(z),M(v))}\;  \overline{R^x_{M(z),M(v)}}.
\end{align}
We also observe 
\begin{equation}
\label{fourr} \Delta_=(u,v)=\sum_{z:\, M(z)=z}R^x_{u,z} \;  q^{\rho(z,v)} \;  \overline{(q-1-x)R^x_{z,M(v)}}= \overline{(q-1-x)} \;  q \;  \Delta_=(u,M(v))= -x \Delta_=(u,M(v))
\end{equation}
and, using the up-down symmetry,
\begin{align}\label{fivve}
 \Delta_=(M(u),M(v))&= \sum_{z:\, M(z)=z}R^x_{M(u),z} q^{\rho(z,M(v))} \;  \overline{R^x_{z,M(v)}}\\
 \nonumber &= (q-1-x)\sum _{z:\, M(z)=z}R^x_{u,z} \;  q^{\rho(z,M(v))} \;  \overline{R^x_{z,M(v)}}\\
 \nonumber &=(q-1-x) \; \Delta_=(u,M(v)).
\end{align}


Eqs. \eqref{onne}, \eqref{twwo},  \eqref{thrree}, \eqref{fourr} and \eqref{fivve} imply 
\begin{align*}
 \Delta(u,v)&=\Delta_{\lhd}(u,v)+\Delta_{\rhd}(u,v)+\Delta_{=}(u,v)\\
 &=-(q-1)\Delta_{\lhd}(u,M(v))-\Delta_{=}(M(u),M(v))+(1-q)\Delta_{\rhd}(u,M(v))- x \Delta_{=}(u,M(v))\\
 &=(1-q)\Delta(u,M(v)),
 \end{align*}
 and the result follows by induction hypothesis since $u\neq M(v)$.

Now suppose $M(u)\rhd u$. We have

\begin{align*}
\Delta_\lhd(u,v)&=  \sum_{z:\, M(z)\lhd z}R^x_{u,z} \; q^{\rho(z,v)} \; \overline{R^x_{z,v}}\\
&=\sum_{z:\, M(z)\lhd z}\big((q-1)R^x_{u,M(z)}+qR^x_{M(u),M(z)}\big) \; q^{\rho(z,v)} \; \overline{R^x_{M(z),M(v)}}\\
&=(q-1)\Delta_\rhd (u,M(v))+q\Delta_\rhd(M(u),M(v)).
\end{align*}

%

Furthermore
\begin{align*}\Delta_\rhd (u,v)&= \sum_{z:\, M(z)\rhd z}R^x_{u,z} \; q^{\rho(z,v)} \; \overline{R^x_{z,v}}\\
&=\sum_{z:\, M(z)\rhd z}R^x_{u,z} \; q^{\rho(z,v)} \; \overline{((q-1)R^x_{z,M(v)} +qR^x_{M(z),M(v)})}\\
&=(1-q) \Delta_{\rhd}(u,M(v))+q\Delta_\lhd (M(u),M(v))
\end{align*}
and (by the up-down symmetry) 
\begin{align*}\Delta_=(u,v)&=\sum_{z:\, M(z)=z}R^x_{u,z} \; q^{\rho(z,v)} \; \overline{R^x_{z,v}}\\
&= q\sum_{z:\, M(z)=z}R^x_{u,z} \;q^{\rho(z,M(v))} \; \overline{(q-1-x) R^x_{z,M(v)}}\\
&=q\sum_{z:\, M(z)=z}R^x_{M(u),z} \; q^{\rho(z,M(v))} \; \overline{R^x_{z,M(v)}}\\
&=q\Delta_=(M(u),M(v)).
\end{align*}
Hence $\Delta(u,v)= q \Delta(M(u),M(v))$, and the result follows  by the induction hypothesis since $M(u)\neq M(v)$.

Finally, we  suppose $M(u)= u$.
We have
\begin{align*}
\Delta_\lhd(u,v)&=\sum_{z:\, M(z)\lhd z}R^x_{u,z} \; q^{\rho{(z,v)}} \;\overline{R^x_{z,v}}\\
&= \sum_{z:\, M(z)\lhd z}(q-1-x)R^x_{u,M(z)} \; q^{\rho(M(z),M(v))} \; \overline{R^x_{M(z),M(v)}}\\
&=(q-1-x)\Delta_{\rhd}(u,M(v))
\end{align*}
and 
\begin{align*}
\Delta_\rhd(u,v)&= \sum_{z:\, M(z)\rhd z}R^x_{u,z} \; q^{\rho(z,v)} \; \overline{R^x_{z,v}}\\
&=\sum_{z:\, M(z)\rhd z } R^x_{u,z} \; q^{\rho(z,v)} \; \overline{((q-1)R^x_{z,M(v)}+qR^x_{M(z),M(v)})}\\
&=(1-q)\Delta_\rhd(u,M(v))+\sum_{z:M(z)\rhd z}\overline{(q-1-x)}R^x_{u,M(z)} \; q^{\rho(M(z),M(v))+2} \; \overline{q R^x_{M(z),M(v)}}\\
&= (1-q)\Delta_\rhd(u,M(v))-x\Delta_{\lhd}(u,M(v))
\end{align*}
and
\begin{align*}
\Delta_=(u,v)&=\sum_{z:\, M(z)=z}R^x_{u,z} \; q^{\rho(z,v)} \; \overline{R^x_{z,v}}\\
&= \sum_{z:\, M(z)=z}R^x_{u,z} \; q^{\rho(z,v)}\; \overline{(q-1-x)R^x_{z,M(v)}}\\
&=q\overline{(q-1-x)}\Delta_=(u,M(v))\\
&=-x\Delta_=(u,M(v)).
\end{align*}
Hence $\Delta(u,v)= -x  \Delta(u,M(v))$, and the result follows by induction hypothesis since $u\neq M(v)$.
\end{proof}

 If the family $\{R^x_{u,v}(q)\}_{u,v\in P}$ of  $R^x$-polynomials defines a $P$-kernel, then we denote the corresponding Kazhdan--Lusztig--Stanley polynomials by $\{P^x_{u,v}(q)\}_{u,v\in P}$. Recall that $P^x_{u,v}(q)$ is zero unless $u\leq v$. 

Note that, by Theorem~\ref{arePkernel}, the Kazhdan--Lusztig--Stanley polynomials $P^x_{u,v}$ are defined if the $R^x$-polynomials satisfy the up-down symmetry.

The following lemma gives a tool to prove the up-down symmetry.

\begin{lem}
	\label{q-1-x}
	Let $(P,S)$ be a  pircon system. Suppose that, for all  $M\in S$ and $z\in P\setminus \{\hat{0}_P\}$ such that $M(z)=z$, there exists  $N\in S$ such that $N(z)\lhd z$ and the restrictions of $M$ and $N$ to $P_{\leq z}$ are strictly coherent. Then the family of Kazhdan--Lusztig $R^x$-polynomials of  $(P,S)$  satisfies the up-down symmetry. 
\end{lem}
\begin{proof}
By Lemma~\ref{cesena}, we need to show that, for all   $M\in S$ and $z\in P$ such that $M(z)=z$, the equation 
\[
R_{u,z}^{x}=(q-1-x) \, R_{M(u),z}^{x}
\]
holds	for all $u\in P$ with $u\leq z$ and $M(u)\rhd u$.

Denote by \( \mathfrak{M}^{M} \) the free ${\mathbb Z}[q]$-module  with $\{m_v : M(v) \text{ is defined} \}$ as a basis: $\mathfrak{M}^{M}= \bigoplus _{v }{\mathbb Z}[q]m_v.$
	 Consider the endomorphism $T_M$ of $\mathfrak M^{M}$ given by
		\[T_M (m_v)= \left\{ \begin{array}{ll}
		m_{M(v)}, & \mbox{if $M(v)  \lhd v$,} \\
		q\,m_{M(v)}+(q-1)\,m_v, & \mbox{if $M(v) \rhd v$,} \\
		(q-1-x)\,m_v, & \mbox{if $M(v) = v$,} 
		\end{array} \right. \] 
	
		Given $f=\sum_{v} f_v(q)\,m_v \in \mathfrak{M}^{M}$ and $w\in P$, we denote   the polynomial
	\( \sum_{v} f_v(q)R^x_{v,w}(q) \in \mathbb Z[q]\) by $f^w$. In this notation, the property of $M$ being strongly calculating reads
	\begin{equation*}
	m_v^w= (T_M (m_v))^{M(w)}
	\end{equation*}
	for all  $v,w\in P$ such that $M(v)$ and $M(w)$ are defined, and $M(w)\lhd w$. It follows, that for all $f\in \mathfrak{M}^{M}$ and all $w$ such that $M(w)\lhd w$, the equation
	
	\begin{equation}\label{recurf}
	f^w=(T_M(f))^{M(w)}
	\end{equation}
holds. 	Since $\frac{q}{-x}=q-1-x$, we need to prove 
	\begin{equation*}\label{enunc}
	(T_M(m_u))^z=(q-1-x)m_u^z
	\end{equation*}

for all $u\leq z$ such that $M(u)\rhd u$. We actually show the (a priori stronger but indeed equivalent) property
	\begin{equation}\label{enuncf}
	(T_M(f))^z=(q-1-x)f^z
	\end{equation}
for all $f\in \mathfrak M ^{M}$; let us proceed  by induction on $\rho(z)$.

Since $(T_M(m_v))^z=(q-1-x)m_v^z=0$ whenever $v\not\leq z$, we may suppose $f\in \mathfrak M ^{\leq z}= \bigoplus _{v\in P_{\leq z}}{\mathbb Z}[q]m_v$. Fix $N\in S$ such that $N(z)\lhd z$ and the restrictions of $M$ and $N$ to $P_{\leq z}$ are strictly coherent.

	Let   $\mathcal O$ be the orbit  of $z$ under the action of $\langle M,N \rangle$, which is a
chain of rank $m(\mathcal O)-1$ by Lemma~\ref{noorbita}. 	

	By \cite[Theorem~4.3]{Mpirc}, the following braid relation 
		\begin{equation}\label{braidrel}\underbrace{\cdots T_M  T_N  T_M}_{m(\mathcal O)}(f)= \underbrace{\cdots T_N  T_M  T_N}_{m(\mathcal O)}(f)
		\end{equation}
holds	for all $f\in \mathfrak M ^{\leq z}$.
	
If the bottom element of $\mathcal O$ is the bottom element of $P$, then $P_{\leq z}$ is a chain by Lemma~\ref{noorbita} and the result follows by Example~\ref{R-chain}. We can therefore suppose that the bottom elements of $\mathcal O$ and of $P$ do not coincide. 
	
	Since   $N(z)\lhd z$, $MN(z)\lhd N(z)$, $NMN(z)\lhd MN(z)$ and so on until we reach the bottom element $\underbrace{\ldots NMN}_{m(\mathcal O)-1}(z)$ of the chain-like orbit $\mathcal O$, we have 
	$$(T_M (f))^z=(T_N T_M (f))^{N(z)} =(T_M  T_N  T_M (f))^{MN(z)}= \cdots =(\underbrace{\cdots T_M  T_N  T_M}_{m(\mathcal O)} (f))^{\underbrace{_{\cdots NMN}}_{m(\mathcal O)-1}(z)}$$
	by Eq. \eqref{recurf} applied repeatedly to both $M$ and $N$.

	Let 
	$$L=   \left\{ \begin{array}{ll}
	M, & \mbox{if $m(\mathcal O)$ is even,} \\
	N, & \mbox{if $m(\mathcal O)$ is odd.}
	\end{array} \right. $$
	The element $\underbrace{\cdots NMN}_{m(\mathcal O)-1}(z)$ is fixed by $L$ and $\underbrace{\cdots T_N  T_M  T_N}_{m(\mathcal O)-1} (f)\in \mathfrak M ^{L}$. By  the induction hypothesis, we have
	\begin{eqnarray*}
		(\underbrace{\cdots T_M  T_N  T_M}_{m(\mathcal O)} (f))^{\underbrace{_{\cdots NMN}}_{m(\mathcal O)-1}(z)}&=&(T_L(\underbrace{\cdots T_N  T_M  T_N}_{m(\mathcal O)-1} (f)))^{\underbrace{_{\cdots NMN}}_{m( \mathcal O)-1}(z)}\\
		&=&(q-1-x) \, (\underbrace{\cdots T_N  T_M  T_N}_{m(\mathcal O)-1}  (f))^{\underbrace{_{\cdots NMN}}_{m(\mathcal O)-1}(z)}\\
		&=& (q-1-x) \,f^z
	\end{eqnarray*}
	where the last equality follows by Eq.~\eqref{recurf}. Hence Eq.~\eqref{enuncf}  follows.
\end{proof}

 We show next that Lemma~\ref{q-1-x} implies that two interesting examples of $R^x$ polynomials such as the Kazhdan--Lusztig--Vogan polynomials associated with the action of  $\Sp(2n,\mathbb C)$  on the flag variety of $\SL(2n,\mathbb C)$ and the parabolic Kazhdan--Lusztig polynomials introduced by Deodhar \cite{Deo87} satisfy the up-down symmetry.
 
The orbits of the action of  $\Sp(2n,\mathbb C)$  on the flag variety of $\SL(2n,\mathbb C)$ are parametrized by the set  
$\iota$ of the \emph{twisted identities} of the symmetric group $S_{2n}$, which is the set $\{\theta (w^{-1} ) w : w\in S_{2n}\}$, where $\theta$ is the involutive automorphism of $S_{2n}$ sending the transposition $s_i=(i,i+1)$ to the transposition $s_{2n-i}=(2n-i,2n-i+1)$, for all $i\in [1,n]$ (we refer the reader to  \cite{AH} for more details on this subject).
The associated Kazhdan--Lusztig--Vogan $R$-polynomials and $Q$-polynomials are indexed by pairs of elements in $\iota$.
 In \cite[Theorem~6.2]{Mpirc}, it is shown that $\iota$, with the order induced by Bruhat order, is a dircon, and that the Kazhdan--Lusztig--Vogan $R$-polynomials and $Q$-polynomials for the action of  $\Sp(2n,\mathbb C)$  on the flag variety of $\SL(2n,\mathbb C)$ coincide, respectively, with the Kazhdan--Lusztig $R^q$-polynomials and $R^{-1}$-polynomials of $\iota$ as a dircon. Since $\iota$ is a dircon,  all special partial matchings are strongly calculating. The map $u\mapsto u*s_i$ given by $u*s_i=\theta(s_{i})us_i$ is a  special partial matching on $\iota$: following   \cite{AH}, we refer to  a special partial matching of this form as a \emph{conjugation matching}.
 

 \begin{thm}
 \label{smileKLV}
The families of  Kazhdan--Lusztig--Vogan $R$-polynomials $\{R_{u,v}\}_{u,v\in \iota}$ and $Q$-polynomials $\{Q_{u,v}\}_{u,v\in \iota}$ for the action of  $\Sp(2n,\mathbb C)$  on the flag variety of $\SL(2n,\mathbb C)$ satisfy the up-down symmetry.
 \end{thm}
\begin{proof}
We use Lemma~\ref{q-1-x}.  Let $z\in \iota \setminus\{e\}$, and $M$ be a special partial matching  such that $M(z)=z$. Since $M$ has fixed points, \cite[Proposition~4.8]{AH} implies that $M$ is a conjugation matching, say $M(v)=v\ast s_i$ for all $v\in \iota$ with $v\leq w$. Since $M(z)=z$, \cite[Theorem~4.3]{AH} implies $s_i\notin D_R(z)$. Hence, since $z\neq e$, there exists $j$ with $j\neq i$ such that $s_j\in D_R(z)$ and let $N$ be the corresponding conjugation special partial matching of $z$.

Notice that $m(\mathcal O)$ is a divisor of $m(s_i,s_j)$, for any orbit $\mathcal O$ of $\langle M,N\rangle $. 
The orbit $\langle M,N\rangle (z)$ is chain-like and 
$$m(\langle M,N\rangle (z))=   \left\{ \begin{array}{ll}
3, & \mbox{if $|i-j|=1$,} \\
2, & \mbox{if $|i-j|>1$.}
\end{array} \right. $$

The orbit $\langle M,N\rangle (u)$ is either chain-like or dihedral and 
$$m(\langle M,N\rangle (u))=   \left\{ \begin{array}{ll}
\mbox{3 or 1}, & \mbox{if $|i-j|=1$,} \\
\mbox{2 or 1}, & \mbox{if $|i-j|>1$.}
\end{array} \right. $$
Hence, we may conclude by Lemma~\ref{q-1-x}.
\end{proof}

\begin{rem}
The fact that Kazhdan--Lusztig--Vogan $R$-polynomials $\{R_{u,v}\}_{u,v\in \iota}$ and $Q$-polynomials $\{Q_{u,v}\}_{u,v\in \iota}$ for the action of  $\Sp(2n,\mathbb C)$  on the flag variety of $\SL(2n,\mathbb C)$ are $\iota$-kernels is an immediate consequence of Theorems~\ref{arePkernel} and~\ref{smileKLV}. The Kazhdan--Lusztig--Stanley polynomials of the Kazhdan--Lusztig--Vogan $Q$-polynomials are the Kazhdan--Lusztig--Vogan $P$-polynomials (see \cite[Eq.~(2)]{AH}). 
\end{rem}

We fix an arbitrary Coxeter system $(W,S)$. Recall that, given $w\in W$, we say that $M$ is a matching of $w$ if $M$ is a matching of the lower Bruhat interval $[e,w]$.
For each $s\in D_{L}(w)$, we have  a special matching 
$\lambda_{s}$ of $w$ 
defined by    $\lambda_{s}(u)=su$, for all $u \in [e,w]$ (see  \cite{CM1} and \cite{CM2}  for more details concerning special matchings of Coxeter systems). 
We call these matchings \emph{left multiplication
matchings}.

Let $H$ be an arbitrary subset of $S$.
Let $w \in W^H$. 
As in \cite{Mtrans, M}, we say that a special matching of $[e,w]$ is {\em $H$-special} provided that 
$$u \leq w, \; u \in W^H, \; M(u) \lhd u \Rightarrow M(u) \in W^H.$$
Note that a $\emptyset$-special matching is a special matching and that a left multiplication matching  is $H$-special 
for all $H \subseteq S$. 

By \cite[Theorem~7.7]{AHH}, the parabolic quotient $W^H$ is a pircon. Recall from \cite{Mpirc} that 
an $H$-special matching $M$ of $[e,w]$ gives rise to a special partial  matching $M^H$ of $[e,w]^H$, which is defined as follows:
\begin{equation}
\label{ricetta}
 M^H(u)= \left\{ \begin{array}{ll}
M(u), & \mbox{if $M(u)  \in W^H $,} \\
u, & \mbox{if $M(u)  \notin W^H $,} 
\end{array} \right. 
\end{equation}
for all $u\in [e,w]^H$. We call \emph{left multiplication partial matchings} the special partial matchings coming from left multiplication matchings.

Let $\SPM^H$
be the set of all special partial matchings of elements in $W^H$ that are obtained from $H$-special matchings by the recipe of Eq.~\eqref{ricetta}: $$\SPM^H=\{M^H \in \SPM_w : M  \text{ is an $H$-special matching of some $w\in W^H\setminus\{e\}$}\}.$$ 

The following result is needed in the proof of Theorem~\ref{up-down parabolica}.
\begin{lem}
\label{finito}
Let $z\in W^H$ and $s,r\in S$. If $rz\lhd z$, $sz\rhd z$ and $sz \notin W^H$, then $r,s \in D_L(sz)$ and $m(s,r)<\infty$. 
\end{lem}
\begin{proof}
Since $z\in W^H$ and $z\lhd sz \notin W^H$, we have $sz=zh$ for a certain $h\in H$. So $r,s \in D_L(sz)$ which, by well known facts, implies $m(s,r)<\infty$. 
\end{proof}

\begin{thm}
\label{up-down parabolica}
Let $(W,S)$ be any Coxeter system, and $H\subseteq S$. Let $\mathcal S^H$ be the set of left multiplication partial matchings of $W^H$. Then
\begin{itemize}
\item $(W^H, \mathcal S^H)$ is a pircon system,
\item the families of Kazhdan--Lusztig $R^{x}$-polynomials of $(W^H, \mathcal S^H)$  satisfy the up-down symmetry.
\end{itemize}
 \end{thm}
\begin{proof}
The properties of a pircon system in Definition~\ref{sistema} are easy to check except the fourth one: let us prove it. Let $u,w\in W^H$ with $u\leq w$,  and $M,N\in \mathcal S^{H}$ such that $M(w)\lhd w$ and $N(w)\lhd w$. Let $s,r\in S$ be such that $M=\lambda_s^H$ and $N=\lambda_r^H$. The orbit  $\langle M,N\rangle (w)$ is dihedral with $m(\langle M,N\rangle (w))=m(s,r)$, and the orbit  $\langle M,N\rangle (u)$ is dihedral, chain-like or a singleton, with  $m(\langle M,N\rangle (u))\in \{1,m(s,r)\}$  by \cite[Lemma~6.8]{Mpirc}. Hence $m(\langle M,N\rangle (u))$ is a divisor of $m(\langle M,N\rangle (w))$.

In order to prove the second statement, we apply Lemmas~\ref{q-1-x}.


Let $s\in S$ and $z \in W^H\setminus\{e\}$ such that   $\lambda_s^H(z)=z$. 
Since $\lambda_s^H(z)=z$, clearly $s\notin D_L(z)$. Hence, since $z\neq e$, there exists $r\in S$ with $r\neq s$ such that $r\in D_L(z)$.  Lemma~\ref{finito} implies $m(s,r)<\infty$.

 The orbit $\langle \lambda_s^H,\lambda_r^H\rangle (z)$ is chain-like with  
$m(\langle \lambda_s^H,\lambda_r^H \rangle (z))= m(s,r)$ by \cite[Lemma~6.8]{Mpirc}.  The orbit $\langle \lambda_s^H,\lambda_r^H \rangle (u)$ is either chain-like or dihedral and, in both cases,  
$m(\langle \lambda_s^H,\lambda_r^H \rangle (u))= m(s,r)$.
Hence, we may conclude by Lemma~\ref{q-1-x}.

\end{proof}

The following result was shown by Brenti  \cite[Corollary~4.2]{BJoA2} for finite Coxeter groups with a complete different proof (which is based on the existence of a longest element in finite Coxeter groups). We prove it for all  Coxeter groups as a direct consequence of the results in this section.

Let $\mathcal M$ be any refinement of the pircon $W^H$ with matchings in $\SPM^H$. By  \cite[Theorem~1.5]{M}, Deodhar's parabolic Kazhdan--Lusztig $R^{x,H}$-polynomials $\{R_{u,v}^{x,H}\}_{u,v\in W^H}$ coincide with the Kazhdan--Lusztig $R^{x,H}$-polynomials  of the refined pircon $(W^H,\mathcal M)$, and, in particular,  with the Kazhdan--Lusztig $R^{x}$-polynomials  of the pircon system $(W^H,\mathcal S^H)$, where  $\mathcal S^H$ is the set of left multiplication partial matchings of $W^H$.

\begin{cor}
\label{brenti}
Let $(W,S)$ be any Coxeter system, and $H\subseteq S$. The families of Deodhar's parabolic Kazhdan--Lusztig $R^{x,H}$-polynomials $\{R_{u,v}^{x,H}\}_{u,v\in W^H}$ satisfy 
	\begin{equation*}
	R_{u,w}^{x,H}=(q-1-x) \, R_{su,w}^{x,H}
	\end{equation*}
	for all $s\in S$ and all $u,w\in W^H$ such that $u\lhd su\in W^H$ and $w\lhd sw\notin W^H$.
 \end{cor}

\section{Deodhar's Duality (revisited)}
\label{deodhar rivisitato}
The purpose of this section is to revisit the results in \cite{Deo91} where Deodhar studies the relationship between the two Hecke algebra modules introduced in \cite{Deo87} and used to define the  parabolic Kazhdan--Lusztig $R^q$-polynomials and $R^{-1}$-polynomials. Not only we generalize Deodhar's results to the setting of pircons, but, using the up-down symmetry, we also shed new light on the case of parabolic Kazhdan--Lusztig polynomials and, in particular, we find the missing involution that was desired in Deodhar's approach (see \cite[Remark~2.5]{Deo91}).

Throughout this section and the next one, we fix a pircon system $(P,S)$ with the property that $S$ consists of quasi special partial matchings of the whole pircon $P$, and whose $R^x$-polynomials satisfy the up-down symmetry.  We consider the free $\mathbb Z[q^{1/2},q^{-1/2}]$-module $\mathcal M_P$ given by
\[
\mathcal M_P=\bigoplus_{u\in P}\mathbb Z[q^{1/2},q^{-1/2}]m_u.\]

Let $(W_{P},S)$ be the Coxeter system given by $m(M,N)=\min ( \{k>0:\, (MN)^k(u)=u\,\,\textrm{ for all }u\in P\}  \cup \{+\infty\})$, for all $M,N\in S$ with $M\neq N$ (i.e., $m(M,N)$ is the order of $MN$ as a permutation of $P$). We denote by $\mathcal H_P$ the Hecke algebra of $W_P$. 

The action of $W_P$ on $P$ can be extended  to an action of the Hecke algebra $\mathcal H_P$ on $\mathcal M_P$ in the two following ways. 
\begin{thm}\label{T_M:action}
Let $x\in \{q,-1\}$. The maps \[
	T_M   \x  m_u= \left\{ \begin{array}{ll}
	m_{M(u)}, & \mbox{if $M(u)  \rhd u$,} \\
	q\,m_{M(u)}+(q-1)\,m_u, & \mbox{if $M(u) \lhd u$,} \\
	x\,m_u, & \mbox{if $M(u) = u$,} 
	\end{array} \right. \]
for all $M\in S$, provide a structure of $\mathcal H_P$-module on $\mathcal M_P$.	
\end{thm}
\begin{proof}
Since the relations of the Hecke algebra $\mathcal H_P$ involve at most two special partial matchings, it is enough to prove the thesis for the set $S$ of cardinality at most 2. Hence the proof is similar to  that of \cite[Theorem~4.3]{M}. In the proof of \cite[Theorem~4.3]{Mpirc}, the hypothesis that is used is that the cardinality of every orbit  of $\langle M,N\rangle$ divides $m(M,N)$: here this holds by the definition of $m(M,N)$.
\end{proof}

As analogues of the map $\iota$ of an Hecke algebra (see subsection~\ref{gruppi di coxeter}), we have the following two maps of $\mathcal M_P$.
\begin{defi} \label{def:iota}We	let $\iota^x$ be the $\mathbb Z$-linear endomorphism of $\mathcal M_P$ given by
	\[
	\iota^x(m_v)=\overline q^{\,\rho(v)}\sum_{u\in P}(-1)^{\rho(u,v)}R^x_{u,v}\,  m_u
	\]
	and \[\iota^x(\sum_{v\in P}f_v\,  m_v)=\sum_{v\in P} \overline{f_v} \, \iota^x(m_v).\]
\end{defi}
Notice that $\iota^x(m_v)=\overline q^{\,\rho(v)}\sum_{u: u\leq v}(-1)^{\rho(u,v)}R^x_{u,v}\,  m_u$ holds, since $R^x_{u,v}=0$ unless $u\leq v$.
\begin{pro}
	The map $\iota^x$ is an involution.
	\end{pro}
\begin{proof}
	Indeed\begin{align*}
	\iota^x(\iota^x(m_v))&=\iota^x(\overline q^{\,\rho(v)}\sum_{z\in P}(-1)^{\rho(z,v)}R^x_{z,v}\,m_z)\\
	&= q^{\rho(v)}\sum_{z\in P}(-1)^{\rho(z,v)}\overline{R^x_{z,v}}\,\overline q^{\,\rho(z)}\sum_{u\in P}(-1)^{\rho(u,z)}R^x_{u,z}\,m_u\\
	&=\sum_{u\in P}\big(\sum_{z\in P} R^x_{u,z}q^{\rho(z,v)}\overline{R^x_{z,v}})\big)(-1)^{\rho(u,v)}m_u,
	\end{align*}
	and the assertion follows by Theorems~\ref{arePkernel} and~\ref{Pkernelequiv} since the $R^x$-polynomials satisfy the up-down symmetry.
\end{proof}
The crucial property that we want to show is the following.
\begin{thm}\label{udiota} For all $h\in \mathcal H_P$ and $m\in \mathcal M_P$, we have
\[	\iota^x(h\x m)= \iota (h)\x (\iota^x(m)).\]
\end{thm}
\begin{proof}
	It is enough to show the statement for $h=T_M$, with $M\in S$, and $m=m_v$. We split the proof into three cases depending on the action of $M$ on $v$.
For notational convenience, we let $\varepsilon_{u,v}=(-1)^{\rho(u,v)}$ for all $u,v\in P$.
	
 If  $M(v)\lhd v$ then, by Theorem~\ref{T_M:action}, Definition~\ref{def:iota} and Eq.~\eqref{klperpirconi}, 
 \begin{eqnarray*}
	q^{\rho(v)}\iota^x(T_M\x m_v)&=&q^{\rho(v)}\iota^x \big( qm_{M(v)}+(q-1)m_v \big)\\
	&=& -   \sum_{u\in P}\varepsilon_{u,v} R^x_{u,M(v)}{}m_u + (\overline q-1)\sum_{u\in P}\varepsilon_{u,v}R^x_{u,v}m_u \\
	&=& - \sum_{u:M(u)\rhd u}\varepsilon_{u,v} R^x_{M(u),v}{}m_u - \sum_{u:M(u)\lhd u}\varepsilon_{u,v} \overline q \big(R^x_{M(u),v}{}-(q-1) R^x_{u,v}{} \big)  m_u\\
	&& - \sum_{u:M(u)= u}\varepsilon_{u,v} R^x_{u,M(v)}{}m_u  + (\overline q-1)\sum_{u\in P}\varepsilon_{u,v}R^x_{u,v}m_u.
	\end{eqnarray*}
Therefore
\begin{align*}	
q^{\rho(v)}\iota^x(T_M\x m_v) &=  \!\!\sum_{u:M(u)=u}\!\!\varepsilon_{u,v} (- R^x_{u,M(v)}{} + (\overline q-1)R^x_{u,v}) m_u+  \!\!\sum_{u:M(u)\rhd u}\!\!\varepsilon_{u,v} \big( - R^x_{M(u),v}{} + (\overline q-1) R^x_{u,v} \big)m_u  \\ 
& \hspace{3mm}+  \!\!\sum_{u:M(u)\lhd u}\!\!\varepsilon_{u,v} \big(- \overline q R^x_{M(u),v}{}-(\overline q-1) R^x_{u,v}{} +(\overline q-1) R^x_{u,v} \big)  m_u  \\
&= \!\!  \sum_{u:M(u)=u}\!\!\varepsilon_{u,v} \big(- \overline{(q-1-x)} R^x_{u,v}{} + (\overline q-1)R^x_{u,v}\big) m_u \\
&\hspace{3mm}+ \!\!\sum_{u:M(u)\rhd u}\!\!\varepsilon_{u,v} \big( - R^x_{M(u),v}{} + (\overline q-1) R^x_{u,v} \big)m_u- \!\!\sum_{u:M(u)\lhd u}\!\! \varepsilon_{u,v} \overline q R^x_{M(u),v}{}  m_u    \\
&=  \! \!\!\sum_{u:M(u)=u}\!\!\!\varepsilon_{u,v} \overline x R^x_{u,v}{}m_u -\!\!\!\sum_{u:M(u)\rhd u}\!\!\!\varepsilon_{u,v} \big(R^x_{M(u),v}{} + (1-\overline q) R^x_{u,v} \big)m_u -\!\!\!\sum_{u:M(u)\lhd u} \!\!\! \varepsilon_{u,v} \overline q R^x_{M(u),v}{}  m_u.
	\end{align*}

Similarly, if $M(v) \rhd v$ then 
	\begin{eqnarray*}
	q^{\rho(v)}\iota^x(T_M\x m_v)&=&q^{\rho(v)}\iota^x(m_{M(v)})=-   \overline q \sum_{u\in P}\varepsilon_{u,v} R^x_{u,M(v)}\,m_u,
	\end{eqnarray*}

and if $M(v)=v$ then	
\begin{eqnarray*}
q^{\rho(v)}\iota^x(T_M\x m_v)&=&q^{\rho(v)}\iota^x(x \, m_{v})=\overline x  \sum_{u\in P}\varepsilon_{u,v} R^x_{u,v}{} \, m_u.
	\end{eqnarray*}

\bigskip 
On the other hand, independently from the action of $M$ on $v$, we have
	\begin{align*}
	q^{\rho(v)}&\iota(T_M)\x (\iota^x(m_v))=(\overline qT_M-(1-\overline q)) \x \Big(\sum_{u\in P}\varepsilon_{u,v}R^x_{u,v}{} \,m_u\Big)\\
	&=  \overline q \sum_{u:\, M(u)\lhd u}\varepsilon_{u,v}  R^x_{u,v}{}   \big( q \, m_{M(u)}+ (q-1) \, m_u\big) + \overline q\sum_{u: M(u)\rhd u}\varepsilon_{u,v}R^x_{u,v}{} \, m_{M(u)} \\
	&\hspace{3mm}+  \overline q \sum_{u:\,M(u)= u}\varepsilon_{u,v} R^x_{u,v} {} \, x \, m_{u} - (1-\overline q) \sum_{u\in P}\varepsilon_{u,v}R^x_{u,v}{} \,m_u\\
	&=  \overline q \sum_{u:\,M(u)\rhd u} -\varepsilon_{u,v} \, qR^x_{M(u),v}{}    \, m_{u} + \overline q \sum_{u:\, M(u)\lhd u}\varepsilon_{u,v} (q-1) R^x_{u,v}{}  \, m_u\\
	  &\hspace{3mm}+ \overline q\sum_{u: M(u)\lhd u}-\varepsilon_{u,v}R^x_{M(u),v}{} \, m_{u} +  \overline q \sum_{u:\, M(u)=u}\varepsilon_{u,v} R^x_{u,v} {} \, x \, m_{u} - (1-\overline q) \sum_{u\in P}\varepsilon_{u,v}R^x_{u,v}{} \,m_u\\
&= \sum_{u:\, M(u)\rhd u} -\varepsilon_{u,v} \big (  R^x_{M(u),v}{}   \, +  (1-\overline q) R^x_{u,v}{} \big) \,m_u\\
  & \hspace{3mm}	+  \sum_{u:\, M(u)\lhd u}\varepsilon_{u,v} \big( (1-\overline q) R^x_{u,v}{}  -\overline q R^x_{M(u),v}{}  -  (1-\overline q) R^x_{u,v}{} \big) \, m_u  +   \sum_{u:\, M(u)= u}\varepsilon_{u,v} \big( \overline q  \, x  - 1 +\overline q  \big)  \,R^x_{u,v}{} \,m_u\\
&=  \sum_{u:\, M(u)\rhd u} -\varepsilon_{u,v} \big (  R^x_{M(u),v}{}   \, +  (1-\overline q) R^x_{u,v}{} \big) \,m_u+  \sum_{u:\, M(u)\lhd u} - \varepsilon_{u,v} \overline q R^x_{M(u),v}{}  \, m_u 
+   \sum_{u:\, M(u)= u}\varepsilon_{u,v} \overline x  \,R^x_{u,v}{} \,m_u,
	\end{align*}
	where the last equality follows since $\overline qx-1+\overline q=\overline x$ holds whenever $x\in\{q,-1\}$. In particular $	\iota(T_M)\x (\iota^x(m_v))= \iota^x(T_M\x m_v)$  when $M(v)\lhd v$. 
	
	If $M(v)\rhd v$, then we use recursion Eq.~\eqref{klperpirconi} to obtain 
\begin{align*}	
q^{\rho(v)}\iota(T_M)\x (\iota^x(m_v)) &=  \sum_{u:\, M(u)\rhd u} -\varepsilon_{u,v} \big (  \overline q(R^x_{u,M(v)}-(q-1)R^x_{u,v})   \, +  (1-\overline q) R^x_{u,v}{} \big) \,m_u\\
  & \hspace{3mm}	+  \sum_{u:\, M(u)\lhd u} - \varepsilon_{u,v} \overline q R^x_{u,M(v)}{}  \, m_u 
+   \sum_{u:\, M(u)= u}\varepsilon_{u,v} \overline x \overline{(q-1-x)} \,R^x_{u,M(v)}{} \,m_u\\
&=-   \overline q\sum_{u}\varepsilon_{u,v} R^x_{u,M(v)}\,m_u,
\end{align*}	
	where the last equality follows since $\overline x \overline{(q-1-x)}=-\overline q$. Hence the assertion follows in this case.
	
	If $M(v)=v$, then we use recursion Eq.~\eqref{klperpirconi} and the up-down symmetry to obtain 
	\begin{align*}	
	q^{\rho(v)}\iota(T_M)\x (\iota^x(m_v)) &=  \sum_{u:\, M(u)\rhd u} -\varepsilon_{u,v} \big (  \overline{(q-1-x)} R^x_{u,v}{}   \, +  (1-\overline q) R^x_{u,v}{} \big) \,m_u\\
& \hspace{3mm}	+  \sum_{u:\, M(u)\lhd u} - \varepsilon_{u,v} \overline q (q-1-x) R^x_{u,v}{}  \, m_u 
+   \sum_{u:\, M(u)= u}\varepsilon_{u,v} \overline x  \,R^x_{u,v}{} \,m_u\\
&=\overline x \sum_{u}\varepsilon_{u,v} R^x_{u,v}{} \, m_u.
\end{align*}	
Hence the assertion follows also in this last case.
\end{proof}

Since the family of $R^x$-polynomials satisfy the up-down symmetry, we can define the $P^x$-polynomials by Theorem~\ref{arePkernel}, and consider the following definition, which generalizes the one of the Kazhdan--Lusztig elements $\{C_w\}_{w\in W}$ and $\{C'_w\}_{w\in W}$ of the Hecke algebra of a Coxeter group. From now on, we let $x\in \{q,-1\}$ and $z$ be such that $\{x,z\}=\{q,-1\}$.
\begin{defi} We let
	\[C^x_w=q^{\frac{\rho(w)}{2}}\sum_{v\in P}(-1)^{\rho(v,w)} \; \overline{q}^{\,\rho(v)} \; \overline{P^x_{v,w}} \; m_v
	\]
	and
	\[C'^x_w= \overline q^{\,\frac{\rho(w)}{2}}\sum_{v\in P}P^z_{v,w}\;m_v.
	\]
\end{defi}
Note that the element $C'^x_w$ is defined using $z$-Kazhdan--Lusztig polynomials. By triangularity, both $\{C^x_w\}_{w\in P}$ and $\{C'^x_w\}_{w\in P}$ are  bases of $\mathcal M_P$. Following what is customary in the classical setting of the Hecke algebras, we call  the $C^x_w$ and the $C'^x_w$ the \emph{Kazhdan--Lusztig elements}.

In \cite{K-L}, Kazhdan and Lusztig define an involution $j_{\mathcal H}$ on the Hecke algebra $\mathcal H$ of any Coxeter group $W$ in the following way:

\[
j_{\mathcal H}(aT_w)= \overline a \,  (-\overline{q})^{\rho(w)} \; T_w
\]
for all $w\in W$ and $a\in \mathbb Z[q^{\frac{1}{2}},q^{-\frac{1}{2}}]$. Kazhdan and Lusztig prove  $j_{\mathcal H} \circ \iota = \iota \circ j_{\mathcal H}$ and $j_{\mathcal H}(C_w)=(-1)^{\ell(w)} C'_w$, for all $w\in W$.
In \cite[Remark~2.5]{Deo91}, Deodhar says that ``one does not have an analogue of $j$ for the parabolic situation since the only possible candidate is not a morphism and does not commute with $\iota$". This is correct, although we show that the right generalization of $j$ to the parabolic setting (and more generally to the setting of pircons) should 
\begin{itemize}
\item not be an equivariant morphism from $(\mathcal M_P,\x)$ to $ (\mathcal M_P,\x)$ but, instead, a $j_{\mathcal H}$-twisted equivariant morphism from $(\mathcal M_P,\x)$ to $ (\mathcal M_P,\z)$ (in the sense of (1) of Theorem~\ref{giustaj}),
\item not commute with $\iota^x$, for a fixed $x$, but instead satisfy  (2) of Theorem~\ref{giustaj}, i.e. $\iota^x$ and $\iota^z$ are conjugated by it.
\end{itemize}
Furthermore, it is worth noticing that Deodhar expected the map $j$ to relate the basis element $C^x_w$ with $C'^x_w$ while, instead, it relates  $C^x_w$ with $C'^z_w$ (see Proposition~\ref{relazione2}).

So we let 
$
j_{P}:\mathcal M_P\rightarrow \mathcal M_P$ be the $\mathbb Z$-linear map given by
\[
j_P(am_w)=\overline a \,  (-\overline{q})^{\rho(w)} \; m_w\]
for all $a\in \mathbb Z[q^{\frac{1}{2}},q^{-\frac{1}{2}}]$ and all $w\in P$. 
Note that  $j_P$ is an involution.

\begin{thm}
\label{giustaj}
The following hold:
\begin{enumerate}
\item	
$$j_P(h\x m)=j_{\mathcal H} (h)\z j_P(m)$$
for all $h\in \mathcal H_P$ and $m\in \mathcal M_P$,
\item $$\iota^x \circ j_P=j_P\circ \iota^z.$$
 \end{enumerate}
\end{thm}
\begin{proof}
In order to prove (1), it is enough to prove the statement for $h=T_M$ and $m=m_v$ by Theorem~\ref{T_M:action}. Since
$$
j_P(T_M\x m_v)=
 \left\{ \begin{array}{ll}
	j_P(m_{M(v)})=(-\overline q)^{\rho(v)+1}  \,  m_{M(v)}, & \mbox{if $M(v)  \rhd v$,} \\
	j_P(q\,m_{M(v)}+(q-1)\,m_v)=   \overline q \, (-\overline q)^{\rho(v)-1}  \,  m_{M(v)} + (\overline q-1)  (-\overline q)^{\rho(v)}  \,  m_{v}, & \mbox{if $M(v) \lhd v$,} \\
	j_P(x\,m_v)=\overline{x} (-\overline q)^{\rho(v)}  \,  m_{v}, & \mbox{if $M(v) = v$,} 
	\end{array} \right. 
$$
and
$$	j_{\mathcal H} (T_M)\z j_P(m_v)= -\overline q \, T_M \z (-\overline q)^{\rho(v)}  \, m_v=  \left\{ \begin{array}{ll}
	(-\overline q)^{\,\rho(v)+1} \; m_{M(v)}, & \mbox{if $M(v)  \rhd v$,} \\
	(-\overline q)^{\rho(v)+1} \; [q\,m_{M(v)}+(q-1)\,m_v], & \mbox{if $M(v) \lhd v$,} \\
	(-\overline q)^{\rho(v)+1} \; z\,m_v, & \mbox{if $M(v) = v$,} 
	\end{array} \right. 
	$$
	the assertion follows since $\overline{x}=-\overline q\,z$.
	
	Let us prove (2).	For all $a\in \mathbb Z[q^{\frac{1}{2}},q^{-\frac{1}{2}}]$ and all $w\in P$, we have
	\begin{align*}
	\iota^x(j_P(a \; m_w))&=\iota^x(\overline a \; (-\overline q)^{\rho(w)}\; m_w)=a \; (-q)^{\rho(w)}\; ({\overline q})^{\rho(w)} \sum_{v\in P} (-1)^{\rho(v,w)}\; R^x_{v,w} \; m_v\\&=a \sum_{v\in P}(-1)^{\rho(v)}\; R^x_{v,w}\; m_v
	\end{align*}
	and
	\begin{align*}
	j_P(\iota^z(a\; m_w))&=j_P\big(\overline a \; ({\overline q})^{\rho(w)}\sum_{v\in P}(-1)^{\rho(v,w)}\; R^z_{v,w}\; m_v\big)=a\; q^{\rho(w)}\sum_{v\in P}(-1)^{\rho(v,w)}\; \overline {R^z_{v,w}}\;  (-\overline q)^{\rho(v)}\; m_v\\
		&=a\sum_{v\in P} (-1)^{\rho(v)} \; R^x_{v,w}\; m_v
	\end{align*}
	where we used the identity $(-q)^{\rho(v,w)} \; \overline {R^z_{v,w}}=R^x_{v,w}$ of Proposition~\ref{relazione}. Hence (2) follows. 
\end{proof}

The following proposition provides some properties of the Kazhdan--Lusztig bases of the pircon system $(P,S)$. Notice that Property (1) is new also for the parabolic Kazhdan--Lusztig bases.

\begin{pro}
\label{relazione2}
The following holds:
	\begin{enumerate}
		
		\item $j_P(C^x_w)=(-1)^{\rho(w)} \; C'^z_w$,
		\item $\iota^x(C'^x_w)=C'^x_w$,
		\item $\iota^x(C^x_w)=C^x_w$.
	\end{enumerate}
\end{pro} 
\begin{proof}
	Let us check (1). We have
	\begin{align*}
	j_P(C^x_w)&=j_P \Big( q^{\frac{\rho(w)}{2}}\sum_{v\in P}(-1)^{\rho(v,w)} \; \overline{q}^{\,\rho(v)} \; \overline{P^x_{v,w}} \; m_v \Big)\\   
	&= {\overline q}^{\,\frac{\rho(w)}{2}}\sum_{v\in P} (-1)^{\rho(v,w)} \; q^{\rho(v)} \; P^x_{v,w}\; (-\overline{q})^{\rho(v)} m_v \\  
	&= \overline q^{\,\frac{\rho(w)}{2}} \sum_{v\in P} (-1)^{\rho(w)} \; P^x_{v,w} \; m_v\\
	&= (-1)^{\rho(w)}\; C'^z_w.
	\end{align*}

	Furthermore
	\begin{align*}
	\iota^x(C'^x_w)&= \iota^x\Big ( \overline{q}^{\,\frac{\rho(w)}{2}} \sum_{v\in P} P^z_{v,w}\; m_v \Big)\\
	&=q^{\frac{\rho(w)}{2}} \sum_{v\in P} \overline{P^z_{v,w}} \sum_{u\in P} \overline q^{\,\rho(v)}(-1)^{\rho(u,v)}\; R^x_{u,v}\;  m_u\\
	&=q^{\frac{\rho(w)}{2}} \sum_{v\in P} \overline{P^z_{v,w}} \sum_{u\in P} {\overline q}^{\,\rho(v)}\; (-1)^{\rho(u,v)}\; (-q)^{\rho(u,v)} \; \overline{R^z_{u,v}} \; m_u\\
	&=q^{\frac{\rho(w)}{2}} \sum_{u\in P} {\overline q}^{\,\rho(u)} \big( \sum_{v\in P} \overline{R^z_{u,v}} \; \overline{P^z_{v,w}}\big) \; m_u\\
	&= q^{\frac{\rho(w)}{2}} \sum_{u\in P} {\overline q}^{\,\rho(u)} \; \overline q ^{\rho(u,w)} \; P^z_{u,w} \; m_u\\
	&=C'^x_w,
	\end{align*}
where the third equality follows by Proposition~\ref{relazione} (3) and the fifth equality follows by the fact that the $P^z$-polynomials are the Kazhdan--Lusztig--Stanley polynomials of the $R^z$-polynomilas (see subsection~\ref{KLS}). This proves  (2). 

Part (2) of Theorem~\ref{giustaj} with (1), (2), and the fact  that $j_P$ is an involution imply (3).
\end{proof}

\section{Characterization and Recursion for Kazhdan--Lusztig elements}

In this section, we briefly present both a characterization  for Kazhdan--Lusztig elements and a recursion for Kazhdan--Lusztig elements and polynomials.

Recall from Section~\ref{deodhar rivisitato} our assumptions on the pircon system $(P,S)$ and in particular that the family of $R^x$-polynomials satisfy the up-down symmetry.
\begin{pro}[Characterization]
\label{caratterizzazione}
 Let $D\in \mathcal M_P$ and $w\in P$ be such that
	\begin{enumerate}
		\item $\iota^x (D)=D$,
		\item $D=\overline q^{\,\frac{\rho(w)}{2}}\sum_{v\in P}Q_{v,w} m_v$, with $Q_{v,w}\in \mathbb Z[q]$, $Q_{w,w}=1$ and $\deg Q_{v,w}<\frac{\rho(v,w)}{2}$.
		\end{enumerate}
		 Then $Q_{v,w}=P^z_{v,w}$ for all $v\leq w$, i.e. $D=C'^x_w$.
\end{pro}
\begin{proof}
Applying $\iota^x$ to both sides of (2),   we obtain
\[
D=q^{\frac{\rho(w)}{2}}\sum_{v\in P} \overline{Q_{v,w}}\,\overline q^{\,\rho(v)}\sum_{u\in P} (-1)^{\rho(u,v)}R^x_{u,v}m_u.
\]
Equating the coefficients of $m_u$ on both sides, by Proposition~\ref{relazione} (3), we obtain
\begin{align*}
Q_{u,w}\overline q^{\,\rho(w)}&=\sum_{v\in P} \overline{Q_{v,w}} \, \overline q^{\,\rho(v)}(-1)^{\rho(u,v)}R^x_{u,v}\\
&= \overline q^{\,\rho(u)} \sum_{v\in P} \overline{R^z_{u,v}}\, \overline{Q_{v,w}}
\end{align*}
 for all $u$ such that $u\leq w$.
The result follows by the definition of the polynomials $P^z_{u,w}$.
\end{proof}

\begin{pro}
\label{rec_el}
Let $v,w \in P$ with $v\leq w$,  $M\in S$ with $M(w)\lhd w$, $v'=\min\{v,M(v)\}$, and $v''=\max \{v,M(v)\}$. Let $C'_M= {\overline q}^{\,\frac{1}{2}}(T_M+1)$ be the Kazhdan--Lusztig element of $\mathcal H_P$ associated with $M$.
	Then 
 $$C'^x_w=C'_M \x C'^x_{M(w)}-\sum_{\substack{u:M(u)\leq u \textrm{, if $x=q$}\\ u:M(u)< u \textrm{, if $x=-1$}  }}\mu(u,M(w))C'^x_u,$$
and
$$P^z_{v,w}=P^z_{v',M(w)}+x_vP^z_{v'',M(w)}-\sum_{\substack{u:M(u)\leq u \textrm{, if $x=q$}\\ u:M(u)< u \textrm{, if $x=-1$}  }} \mu(u,M(w)) q^{\frac{\rho(u,w)}{2}} P^z_{v,u},$$
	where  $$x_v=\begin{cases}
	x& \textrm{if $M(v)=v$}\\q& \textrm{otherwise.} 
	\end{cases}$$
\end{pro}
\begin{proof}
We use the characterization of the elements $C'^x_w$ of Proposition~\ref{caratterizzazione}. For short, let $E_w^x$ denotes the element
$$C'_M\x C'^x_{M(w)}-\sum_{\substack{u:M(u)\leq u \textrm{, if $x=q$}\\ u:M(u)<u \textrm{, if $x=-1$}  }}\mu(u,M(w))C'^x_u.$$
Since $\iota(C'_M)=C'_M$, we have  $\iota^x(E_w^x)=E^x_w$, by Theorem~\ref{udiota} and Proposition~\ref{relazione2}. Let us expand $E^x_w$ as a linear combination of the elements $m_v$.
	Since
	\begin{align*}
	C'_M\x C'^{x}_{M(w)}&= \overline q^{\,\frac{1}{2}}(1+T_M)\x \Big(\overline q^{\,\frac{\rho(w)-1}{2}}\sum_{v\in P}P^{z}_{v,M(w)}m_v\Big)\\
	&=\overline q^{\,\frac{\rho(w)}{2}}\big( \sum_{v\in P}P^{z}_{v,M(w)}m_v + \sum_{v:\, M(v)\rhd v}P^{z}_{v,M(w)}m_{M(v)}\\
	& \hspace{3mm}+\sum_{v:\, M(v)\lhd v}P^{z}_{v,M(w)} (qm_{M(v)}+(q-1)m_v) +\sum_{v:\, M(v)=v}xP^{z}_{v,M(w)} m_v \big),
	\end{align*}
the coefficient of $m_v$ in the expansion of $E^x_w$ is $q^{-\frac{\rho(w)}{2}}$ times the following polynomial: 
 \begin{align*}
 P^z_{v,M(w)}+\chi_{M(v)\lhd v}P^z_{M(v),M(w)}+\chi_{M(v)\rhd v}qP^z_{M(v),M(w)}+\chi_{M(v)\lhd v} (q-1)P^z_{v,M(w)}\\+\chi_{M(v)=v}x P^z_{v,M(w)}-\sum_{\substack{u:M(u)\leq u \textrm{, if $x=q$}\\ u:M(u)<u \textrm{, if $x=-1$}  }} \mu(u,M(w))q^{\frac{\rho(u,w)}{2}}P^z_{v,u}
 \end{align*}
 where  $$\chi_F=\begin{cases}
	1& \textrm{if $F$ is true}\\0& \textrm{if $F$ is false.} 
	\end{cases}$$
 The assertion follows by checking that this polynomial coincides with $$P^z_{v',M(w)}+x_vP^z_{v'',M(w)}-\sum_{\substack{u:M(u)\leq u \textrm{, if $x=q$}\\ u:M(u)<u \textrm{, if $x=-1$}  }} \mu(u,M(w)) q^{\frac{\rho(u,w)}{2}} P^z_{v,u}$$
 and has degree smaller than $ \frac{\rho(v,w)}{2}$.
\end{proof}	


\end{document}